\documentclass[12pt,a4paper,authoryear]{article}

\usepackage[top=3cm,bottom=3cm,left=2.4cm,right=2.4cm]{geometry}
\usepackage{amsmath,amssymb,amsfonts,amsthm}
\usepackage[utf8]{inputenc}
\usepackage[english]{babel}
\usepackage[colorlinks]{hyperref}
\usepackage[table]{xcolor}
\usepackage{comment}
\usepackage{booktabs}
\usepackage[round]{natbib}
\usepackage{tikz}
\usepackage{subcaption}

\theoremstyle{plain}  
\newtheorem{theorem}{Theorem}[section]
\newtheorem{definition}[theorem]{Definition}
\newtheorem{defprop}[theorem]{Definition and Proposition}
\newtheorem{lemma}[theorem]{Lemma}
\newtheorem{proposition}[theorem]{Proposition}
\newtheorem{corollary}[theorem]{Corollary}
\newtheorem{remark}[theorem]{Remark}

\AtBeginDocument{%
	\hypersetup{
		breaklinks=true,
		urlcolor = blue,
  		linkcolor = blue,
  		citecolor = orange,
	}
}

\DeclareMathOperator*{\argmin}{arg\,min}

\def\supp{\mathrm{supp}}

\begin{document}


\newcommand{\C}{\mathbb{C}}
\newcommand{\N}{\mathbb{N}}
\newcommand{\R}{\mathbb{R}}
\newcommand{\tf}{\mathcal{F}}
\newcommand{\s}{\tau}


\title{Time-asymptotic propagation of approximate solutions of Schrödinger equations with both potential and initial condition in Fourier-frequency bands}

\author{Florent Dewez\footnote{INRIA, DataShape team-project, Saclay -- Ile-de-France research center, France. E-mail: \texttt{florent.dewez@inria.fr} or \texttt{florent.dewez@outlook.com}}}

\date{\today}

\maketitle


\begin{abstract}
	In this paper, we consider the Schrödinger equation in one space-dimension with potential and we aim at exhibiting dynamic interaction phenomena produced by the potential. To this end, we focus our attention on the time-asymptotic behaviour of the two first terms of the Dyson-Phillips series, which gives a representation of the solution of the equation according to semigroup theory. The first term is actually the free wave packet while the second term corresponds to the wave packet resulting from a first interac\-tion between the free solution and the potential. In order to follow a method developed in a series of papers and aiming at describing propagation features of wave packets, we suppose that both the potential and the initial datum are in bounded Fourier-frequency bands; in particular a family of potentials satisfying this hypothesis is constructed for illustration. We show then that the two terms are time-asymptotically localised in space-time cones which depend explicitly on the frequency bands. Since the inclination and the width of these cones indicate the time-asymptotic motion and dispersion of the two terms, our approach permits to highlight interaction phenomena produced by the potential.
\end{abstract}

\vspace{0.3cm}



\section{Introduction} \label{sec:intro}

In this paper, we aim at describing the time-asymptotic spatial propagation of an ap\-pro\-xi\-mate solution of the Schrödinger equation with potential on the real line, \emph{i.e.}
\begin{equation} \label{eq:schro-potential0}
	\left\{ \begin{array}{rl}
			& \hspace{-2mm} i \, \partial_t u(t) = -\partial_{xx} u(t) + V(x) u(t) \\ [2mm]
			& \hspace{-2mm} u(0) = u_0
	\end{array} \right. \; ,
\end{equation}
for all $t \geqslant 0$, where the initial datum $u_0$ is an element of $H^1(\R)$ and the potential $V$ belongs to the Sobolev space $W^{1,\infty}(\R)$. Under these hypotheses, we show that the solution can be represented by the Dyson-Phillips series whose truncation at order 2 provides the ap\-pro\-xi\-mate solution we consider. In particular, we prove that the two terms of this truncated series can be written as wave packets. This allows to use the method developed in the series of papers \citet{amhdr2012, amd2017,dewez2018,dewez2020} which provide time-asymptotic expansions of wave packets in space-time cones. This method, which can be interpreted as a mathematical formulation of the group velocity principle, is mainly based on the hypothesis of initial data in bounded Fourier-frequency bands, that is to say initial data whose Fourier transform is supported on a bounded interval. Under the hypothesis that the potential is also in a Fourier-frequency band, we apply this method to each term of the truncated series together with additional arguments to obtain time-asymptotic expansions in cones. These expansions highlight the influence of the frequencies of the potential on the propagation of the approximate solution, exhibiting in particular reflection and advanced or retarded transmission phenomena.\\

To obtain our results, we consider the setting of semigroups \citep{en2000} where the solution of a linear evolution equation is interpreted as the orbit of the initial datum under the action of a semigroup. In the case of evolution equations perturbed by bounded operators, such as the Schrödinger equation \eqref{eq:schro-potential0} described above, classical results from semigroup theory claim that the dynamics can be still described by a semigroup. In particular, this perturbed semigroup is the limit in the operator norm of a series called the Dyson-Phillips series, whose terms are iteratively defined; see Appendix \ref{sec:appendix} for more precise details.

This series has been used in the literature to derive some asymptotic properties of evolution equations. For instance, it has been used in \cite{al2014} to characterise honest solutions of general transport equations; see \citet[Sec. 1.2]{ba2006} for an introduction on honesty theory. As mentioned by the authors, the use of the Dyson-Phillips series in their study is equivalent to an approach based on the resolvent of the underlying operator but appears to be robust enough to be applied to other problems for which the resolvent approach would be inappropriate. The series has also been recently used in \citet{lmk2021} to compute the rate of convergence to e\-qui\-li\-brium of solutions of collisionless kinetic equations with diffuse boundaries. In this setting, each term has a physical meaning: the $n$-th term of the Dyson-Phillips series is actually the solution of the equation after having experienced $n$ rebounds with the boundary.

In the interaction picture in quantum mechanics, this series provides a representation of the unitary time-evolution operator associated with Schrödinger equations, where the free Hamiltonian is perturbed by small perturbations; see \citet[Sec. 3.5]{weinberg1995}. In particular, this permits to represent the scattering matrix (or S-matrix or S-operator), which describes the transition between two states from distant past to distant future, as a series. In particular, the two first terms of the Dyson-Phillips series can be used to derive Fermi's Golden rule in quantum mechanics. This rule provides an approximate but explicit formula for the transition rate from one energy state to another one; see \citet[Sec. 5.6]{sakurai1993}.

It is also noteworthy that, in scattering theory in physics, the idea of truncating a series is also used to study the spatial asymptotic behaviour of the wave resulting from the interaction with a potential: this is the Born approximation \citep[Sec. 7.2]{sakurai1993}, which consists in truncating the series representing the amplitude of the scattered wave. In particular, the $n$-th term of this series corresponds to scattering viewed as a $n$-step process. In the setting of evolution equations on graphs, a similar interpretation has been proposed in \citet{aman2017} where the difference between the absolutely con\-ti\-nuous part of the solution of the Schrödinger equation on a tadpole graph and the solution of the associated Neumann half-line problem is studied. The authors have proved that this difference on the queue of the tadpole is actually given by a series whose $n$-th term can be interpreted as a wave packet passing from the head of the tadpole into its queue after $n$ cycles around the head.\\

Inspired by the above approaches, we consider the truncated Dyson-Phillips series re\-pre\-sen\-ting the time-evolution of the Schrödinger equation \eqref{eq:schro-potential0} and study the spatial pro\-pa\-ga\-tion for large time of the two first terms. More precisely, if $\big(S(t)\big)_{t \geqslant 0}$ denotes the semigroup generated by the operator $\big( i \partial_{xx} - i V(x), H^3(\R) \big)$ then it can be written as the Dyson-Phillips series:
\begin{equation*}
	S(t) = \sum_{n=1}^{+\infty} S_n(t) \; ,
\end{equation*}
which converges here in the operator norm on $H^1(\R)$ and the terms are recursively defined; we refer to Theorem \ref{thm:dyson} for a generic statement on the existence and convergence of this series. Given this, we focus on the following approximate solution of the Schrödinger equation \eqref{eq:schro-potential0}:
\begin{equation*}
	\sum_{n=1}^{2} S_n(t) u_0 = S_1(t) u_0 + S_2(t) u_0 \; .
\end{equation*}
Studying this approximation is relevant not only to present the main lines of our approach in a simple setting and to develop first arguments for the study of the entire series but also to obtain first results which are interpretable from a physical point of view. Using the recurrence relation defining the terms of the series, we prove in this paper that the term $S_n(t) u_0(x)$, for $n=1,2$, can be written as follows
\begin{equation} \label{eq:sn_oi}
	S_n(t) u_0(x) = \frac{1}{2\pi} \int_\R U_n(t,p) \, e^{-itp^2 + ixp} \, dp \; ,
\end{equation}
for all $ (t,x) \in \R \times \R$, where the amplitude $U_n$ can be explicitly computed. As proved in Theorem \ref{thm:existence_unique}, the term $S_1(t) u_0$ is actually the free solution and so its amplitude $U_1$ does not depend on time since it is the Fourier transform of $u_0$. And we refer to Proposition \ref{prop:rewriting} for a formula for $U_2$ which depends explicitly on the potential. Here the second term $S_2(t) u_0$ can be interpreted as the wave packet issued from a first interaction between the free solution and the potential. In particular, the amplitude $U_2$ describes an interaction between the first term of the Dyson-Phillips series, namely the free solution, and the potential through its Fourier transform.\\

The two terms $S_1(t) u_0$ and $S_2(t) u_0$ being given as wave packets, we are in position to exploit the method developed in the papers \citet{amhdr2012, amd2017,dewez2018,dewez2020} to study their time-asymptotic propagation. The authors are interested in describing precisely the propagation of free wave packets\footnote{The expression \emph{free wave packet} refers here to wave packets having time-independent amplitudes as in \eqref{eq:wave_packet}.} of the form
\begin{equation} \label{eq:wave_packet}
	u_f(t,x) = \frac{1}{2 \pi} \int_\R U (p) \, e^{-itf(p) + ixp} \, dp \; ,
\end{equation}
for a sufficiently regular amplitude and compactly supported $U: \R \longrightarrow \C$ and a strictly convex symbol $f: \R \longrightarrow \R$. The method is inspired by the group velocity principle in physics: broadly speaking, this principle claims that the envelope of a free and almost monochromatic wave packet (its Fourier transform is sharply peaked around a certain frequency $\overline{p}$) moves at a speed given by $f'(\overline{p})$, providing then information on the propagation. In particular, the wave packet \eqref{eq:wave_packet} can be interpreted as the solution of the following type of dispersive equation:
\begin{equation} \label{eq:disp_equation}
	\left\{ \begin{array}{l}
			\left[ i \, \partial_t - f \big(D\big) \right] u_f(t) = 0 \\ [2mm]
			u_f(0) = u_0
	\end{array} \right. \; ,
\end{equation}
where $U$ is actually the Fourier transform of the initial datum $u_0$; let us note that the symbol is given by $f(p) = p^2$ in the Schrödinger setting.

To formalise mathematically the group velocity principle, the authors in \citet{amhdr2012} have proposed to consider initial data in frequency bands (and not necessarily localised around a given frequency). To be precise, the Klein-Gordon equation on a star-shaped network is studied in \citet{amhdr2012} and so the notion of frequency in this setting is associated with a Fourier-like transform which diagonalises the operator of the equation; however the ideas of the method remain the same if unitary transforms other than the Fourier transform are used. To obtain then information on the position of the solution for long times, they have approximated it by a spatially localised function which tends to the solution when the time tends to infinity. To obtain such an approximation, they have precisely applied the version of the stationary phase method given in \citet[Thm. 7.7.5]{hormander83} to a wave packet representation of the solution. Since the principle of the stationary phase method consists in approximating wave packets (or more generally oscillatory integrals) by a term including the integrand evaluated at the stationary point $(f')^{-1}\big(\frac{x}{t}\big)$ of the phase, this makes that the approximation is supported in a space-time cone for a compactly supported amplitude. Information on both the motion and the dispersion of the wave packet can then be deduced and visualised via the inclination and the width of the cone.

This approach has been then extended to the case of amplitudes $U$ having integrable singular points in \citet{amd2017}. To do so, the authors use the version of the stationary phase method given in \citet[Sec. 2.9]{erdelyi1956} which covers the case of oscillatory integrals with integrable singular amplitudes. A modern formulation of this version together with a detailed proof is proposed in \citet{amd2017} and the result is applied to the solution of the free Schrödinger equation on the line. It is proved that the existence of an integrable singular Fourier-frequency for the initial datum makes the associated solution time-asymptotically localised around a space-time line given by this frequency. This shows that the solution travels mainly at the speed given by the inclination of the line.

Another tool to study the time-asymptotic propagation of wave packets solutions of equation \eqref{eq:disp_equation} has been developed in \citet{dewez2018}. An extension of the classical van der Corput Lemma \citep[Prop. 2, Chap. VIII]{stein93} to oscillatory integrals with amplitudes having integrable singular points and phases with stationary points of real order is given and used to establish uniform and explicit estimates of the wave packet \eqref{eq:wave_packet} in space-time cones. These results have been motivated by the fact that it is not possible to obtain uniform asymptotic expansions of wave packets in case of singular frequencies as explained in \citet{amd2017}. Hence the precision of an expansion has been removed in favour of more flexible uniform estimates.

Recently a refinement of the approach based on time-asymptotic expansions has been given in \citet{dewez2020}. In this paper, an improvement of the technical arguments leading to expansions in cones of free wave packets is given. This improvement permits to put the origin of the cone at the mean position of the wave packet at the time when its variance is minimal. This result improves three points:
\begin{itemize}
	\item the coefficient of the remainder estimate is minimal for this position of the cone;
	\item the first term of the expansion has the same mean position as the wave packet;
	\item the difference between the variances of the wave packet and the first term is constant.
\end{itemize}
It is worth noting that the original versions of the method proposed in \citet{amhdr2012,amd2017,dewez2018} do not permit such a precision since these versions forced the cone origin to be at the space-time point $(0, 0)$.\\

Our main results in the present paper are given in Corollary \ref{cor:s1_exp} and Theorem \ref{thm:s2_exp}. Corollary \ref{cor:s1_exp} provides a time-asymptotic expansion of the first term $S_1(t) u_0$ of the Dyson-Phillips series in a space-time cone together with an explicit and uniform remainder estimate. An estimate of $S_1(t) u_0$ outside the cone is also provided, showing especially that the time-decay is faster outside the cone than inside. In this result, the Fourier transform of the initial datum $u_0 \in H^3(\R)$ is supposed to be a continuously differentiable function on $\R$ supported on the interval $[p_1, p_2]$, where $p_1 < p_2$ are two real numbers. Technically speaking, Corollary \ref{cor:s1_exp} is a straightforward application of the above mentioned method for time-asymptotic expansions of wave packets; more precisely, we apply Theorem \ref{thm:oi_exp}, which is a slight extension of \citet[Thm. 1.1]{dewez2020}, to the representation \eqref{eq:sn_oi} of $S_1(t) u_0$. \\
Theorem \ref{thm:s2_exp} provides the same kind of result as Corollary \ref{cor:s1_exp} for the second term $S_2(t) u_0$. Here we require in addition the Fourier transform of the initial datum to be a $\mathcal{C}^4$-function and $0 \notin [p_1, p_2]$. Furthermore the Fourier transform of the potential $V \in L^2(\R)$ is supposed to be a $\mathcal{C}^3$-function and to be supported on the interval $[a, b]$, where $a < b$ are two real numbers. It is noteworthy that our approach does not require any self-adjointness argument and so complex-valued potentials can be considered; for a real-valued potential, its Fourier transform is a Hermitian function implying especially that the frequency band is symmetric, \emph{i.e.} $a = -b$. The proof of Theorem \ref{thm:s2_exp} consists in firstly showing that $S_2(t) u_0$ is time-asymptotically close to an explicit free wave packet and secondly in expanding this free wave packet to one term. We decompose then the proof of Theorem \ref{thm:s2_exp} into two steps: first, we show that the amplitude $t \mapsto U_2(t,.)$ converges sufficiently fast to a continuously differentiable function $W_\infty$ as $t$ tends to infinity. The free wave packet whose amplitude is given by $W_\infty$ provides then the desired approximation of $S_2(t) u_0$, the amplitude $W_\infty$ being actually the Fourier transform of a function which may be interpreted as a fictive initial datum. In particular, we prove that the support of $W_\infty$ is contained in $[p_1 + a, p_2 + b]$. Since the amplitude $U_2(t,.)$ describes an interaction between the potential and the free solution, this first step amounts to proving a decreasing interaction over time. In the second step, we apply once again Theorem \ref{thm:oi_exp} but this time to the wave packet whose amplitude is given by $W_\infty$, providing finally a time-asymptotic expansion for $S_2(t) u_0$ in a cone and a uniform and explicit estimate outside. Similar to $S_1(t) u_0$, the time-decay turns out to be slower in the associated cone than outside.

We interpret now from a physical point of view the two preceding results. First of all, as explained before, each of the two terms of the Dyson-Phillips series tends to be time-asymptotically localised in a space-time cone. Corollary \ref{cor:s1_exp} shows that $S_1(t) u_0$ propagates in the cone delimited by the two lines $\frac{x}{t} = 2 p_1$ and $\frac{x}{t} = 2 p_2$. On the other hand, the second term $S_2(t) u_0$ is rather localised in the cone delimited by $\frac{x}{t} = 2 (p_1 + a)$ and $\frac{x}{t} = 2 (p_2 + b)$. Note that the factor $2$ comes from the fact we consider the operator $-\partial_{xx}$ and not $-\frac{1}{2} \partial_{xx}$ as it is classically done in quantum mechanics.\\
To illustrate these results, we refer to Figures \ref{fig:ret_trans} and \ref{fig:refl}. There we give a representation of the cones associated with $S_1(t) u_0$ and $S_2(t) u_0$: the first one is given by the dark grey-shaded cone with solid lines while the second one is delimited by the dashed lines. To understand the influence of the potential on the propagation, we compare the cone associated with the free wave packet $S_1(t) u_0$, which is by definition not influenced by the potential, and the cone associated with $S_2(t) u_0$, which is the wave packet issued from a first interaction with the potential. In both figures, the potential is in a symmetric frequency\footnote{In favour of readability, the expression \emph{frequency} will refer to \emph{Fourier-frequency} throughout the rest of the present paper.} band (namely a band of the form $[-b, b]$) as in the case of real-valued potentials. In this case, it is noteworthy that the cone associated with $S_2(t) u_0$ contains always the cone associated with $S_1(t) u_0$.\\
In the two figures, the initial datum is supposed to be in a positive frequency band but these frequencies are assumed to be larger in Figure \ref{fig:ret_trans} than in Figure \ref{fig:refl}. This implies that the free wave packet in the first case moves faster to the right in space than the free wave packet in the second case. This can be visualised via the cones: the cone with solid lines is more inclined to the right in Figure \ref{fig:ret_trans} than in Figure \ref{fig:refl}.\\
Regarding the cone associated with $S_2(t) u_0$, we observe that the right dashed line is always more inclined to the right than the right solid line. This means that the second term may travel faster than the free wave packet, indicating an advanced transmission phenomenon. On the other hand, the left dashed line is either inclined to the right or to the left:
\begin{itemize}
	\item in the first case, the second term still travels to the right but a retarded transmission may occur since the left solid line is more inclined to the right than the left dashed line;
	\item in the second case, a part of $S_2(t) u_0$ may travel to the left in space, indicating that a part of $S_2(t) u_0$ is reflected by the potential.
\end{itemize}

\begin{figure}
	\centering
	\begin{subfigure}{.475\textwidth}
		\centering
		\begin{tikzpicture}
			\draw[->] (-1,0) -- (5,0) node[below right] {$x$};
			\draw[->] (0,-1) -- (0,5) node[above left] {$t$};
			\clip (-1,-1) rectangle (5,5);
			\draw [fill=gray,fill opacity=0.2] plot [-,samples=100,domain=0:6](\x,{1/1.8*\x}) -- plot [-,samples=100,domain=6:0](\x,{1/1.2*\x}) ;
			\draw [dashed,fill=gray,fill opacity=0.2] plot [samples=100,domain=0:6](\x,{1/2.1*\x}) -- plot [samples=100,domain=6:0](\x,{1/.9*\x}) ;
		\end{tikzpicture}
		\caption{High and positive frequency initial datum.}
		\label{fig:ret_trans}
	\end{subfigure}
	\hfill
	\begin{subfigure}{.475\textwidth}
		\centering
		\begin{tikzpicture}
			\draw[->] (-1,0) -- (5,0) node[below right] {$x$};
			\draw[->] (0,-1) -- (0,5) node[above left] {$t$};
			\clip (-1,-1) rectangle (5,5);
			\draw [fill=gray,fill opacity=0.2] plot [-,samples=100,domain=0:6](\x,{1/.8*\x}) -- plot [-,samples=100,domain=6:0](\x,{1/.2*\x}) ;
			\draw [dashed,fill=gray,fill opacity=0.2] plot [samples=100,domain=0:6](\x,{1/1.1*\x}) -- plot [samples=100,domain=-1:0](\x,{-1/.1*\x}) ;
		\end{tikzpicture}
		\caption{Low and positive frequency initial datum.}
		\label{fig:refl}
	\end{subfigure}
	 \caption{Illustrations of space-time cones associated with the terms $S_1(t) u_0$ and $S_2(t) u_0$ -- The cones associated with $S_1(t) u_0$ and $S_2(t) u_0$ are respectively delimited by solid and dashed lines and the cone origins are put at $(0, 0)$ for the sake of clarity.}
\end{figure}

Theorem \ref{thm:s2_exp} together with Corollary \ref{cor:s1_exp} indicates that, in the case of an initial wave packet moving to the right in space, the positive frequencies of the potential tend to accelerate the motion while the negative ones tend to slow it down, even reverse it. Hence the existence of a retarded transmission or a reflection depends on the frequency bands of both the initial datum and the potential: for an initial datum in the positive frequency band $[p_1, p_2]$ and a potential in the band $[a, b]$, the left dashed line is inclined to the left if and only if $p_1 + a < 0$. Hence an initial datum with high momentum is less likely to produce a reflection than an initial datum with small momentum. With very similar arguments, it is straightforward to extend the above interpretation to the case of negative frequency bands for the initial datum.\\
To finish with the interpretation, we emphasise that the cones share the same origins in the figures for a clearer explanation of the propagation of the two terms. Nevertheless these origins are not necessarily equal to $(0,0)$. Indeed, as explained above, the origin of the cone associated with the free wave packet $S_1(t) u_0$ is put at the mean position of the wave packet at the time when its variance is minimal. Regarding the second term, it is localised in a cone whose origin is given by the mean position of the free wave packet whose amplitude is $W_\infty$ at the time when the variance is minimal. Since $W_\infty$ can be interpreted as the Fourier transform of a fictive initial datum, it would be an interesting issue to make explicit the dependence of the second cone origin on the potential $V$ and on the true initial datum $u_0$ to understand better this fictive initial datum.

It is noteworthy that our results show that, under our assumptions on $u_0$ and $V$, the second term $S_2(t) u_0$ behaves as a free wave packet for long times. One reason for this pro\-per\-ty may be the following: we recall that Theorem \ref{thm:s2_exp} holds for initial data in a frequency band $[p_1, p_2]$ such that $0 \notin [p_1, p_2]$, together with additional regularity assumptions on $\widehat{u_0}$ and $\widehat{V}$ (that is to say spatial decay assumptions on $u_0$ and $V$). From a physical point of view, this means that the initial wave packet has no zero-velocity component: the pro\-ba\-bi\-li\-ty of finding the associated free wave packet in a given bounded spatial interval tends to 0 as $t$ goes to infinity. One may expect that the absence of zero-velocity components of a wave packet interacting with a short range potential may help to prevent trapped components after a first interaction.\\

Since a natural extension of the present work is to apply our methodology to all the terms of the series, an important step will be to determine sufficient hypotheses on both the initial datum and the potential to assure that each term of the series behaves time-asymptotically as a free wave packet. This could lead to information on the propagation of the true solution. From a technical point of view, the two main challenging tasks of this extension will be firstly to rewrite each term as a wave packet of the form \eqref{eq:sn_oi} with explicit amplitude and secondly to study the convergence of the series of expansions and the series of remainders. A thorough comparison of these future results based on the Dyson-Phillips series with the classical literature based on spectral analysis, such as \citet{weder2000,gs2004,ekmt2016}, would be constructive.

In the present paper, we have focused our attention on potentials in frequency bands to make explicit the influence of its frequencies on the propagation of the second term of the Dyson-Phillips series. From a technical point of view, this has been achieved by expanding time-asymptotically this second term and by proving that the resulting expansion lies in a space-time cone whose inclination depends explicitly on the potential frequencies. Another extension of our approach will be then to consider potentials which are no longer in frequency bands. Nevertheless, in this setting, the second term of the Dyson-Phillips series is no longer in a frequency band, preventing its time-asymptotic expansion from being localised in a cone. Hence another tool to describe the dynamics of a wave packet which is not in a frequency band would be required. In view of this, one may use weighted $L^p$-norms which provide some information on the spatial localisation. These norms have already been exploited to establish dispersive estimates for the Schrödinger equation \eqref{eq:schro-potential0} under some hypotheses on the decay of the potential; we refer for instance to \citet{goldberg2007,schlag2007,eht2015}. For our purpose, it would be relevant to make time-dependent the weight to take into account the motion and the dispersion of the wave packet.\\


The paper is organised as follows: the following section is devoted to the well-posedness in $H^1(\R)$ of the Schrödinger equation with potential \eqref{eq:schro-potential0} for initial data in $H^3(\R)$ and potentials in $W^{1,\infty}(\R)$, and to the Dyson-Phillips series representing the solution. In Section \ref{sec:potential}, we introduce the frequency band hypothesis for the potentials and provide an explicit family of such potentials. In Section \ref{sec:s1}, we focus on the time-asymptotic behaviour of the first term of the Dyson-Phillips series; to do so, we give a slight extension of \citet[Thm. 1.1]{dewez2020} which is then applied to the first term to obtain a time-asymptotic expansion in a cone and a uniform estimate outside. The time-asymptotic behaviour of the second term of the Dyson-Phillips series is studied in Section \ref{sec:s2}; we first provide a representation of it as a wave packet with time-dependent amplitude, we show then that this amplitude tends to a constant function when the time tends to infinity and we finally apply once again the extension of \citet[Thm. 1.1]{dewez2020} to the free wave packet with the limit amplitude to get the result on the second term. And Appendix \ref{sec:appendix} provides some classical but useful results from classical and functional analyses.


\section{Dyson-Phillips series for the Schrödinger equation with potential} \label{sec:dyson}

Let us recall the Schrödinger equation with potential on the line, \emph{i.e.}
\begin{equation} \label{eq:schro-potential2}
	\left\{ \begin{array}{rl}
			& \hspace{-2mm} i \, \partial_t u(t) = -\partial_{xx} u(t) + V(x) u(t) \\ [2mm]
			& \hspace{-2mm} u(0) = u_0
	\end{array} \right. \; ,
\end{equation}
for all $t \geqslant 0$, where the potential $V$ is an element of the Sobolev space $W^{1, \infty}(\R)$; we emphasise that the potential is not required to be real-valued in this paper. The aim of the present section is to assure existence and uniqueness of a solution for the Schrödinger equation \eqref{eq:schro-potential2} in $H^1(\R)$ by exploiting the theory of semigroups. We also introduce the Dyson-Phillips series which provides a representation of the solution as a series.\\

Before stating these results, let us define some objects that will be used throughout the rest of this paper.

\begin{definition} \label{def:misc}
	\begin{enumerate}
		\item Let $\tf_{x \rightarrow p}$ denote the Fourier transform on $L^2(\R)$ and $\tf_{p \rightarrow x}^{-1}$ its inverse. For $f$ in the Schwartz space $\mathcal{S}(\R)$, the element $\tf_{x \rightarrow p} f$ defines a complex-valued function on $\R$ given by
		\begin{equation*}
			\forall \, p \in \R \qquad \big(\tf_{x \rightarrow p} f\big) (p) = \int_{\R} f(x) \, e^{-ixp} \, dx \; .
		\end{equation*}
		If there is no risk of confusion, we shall note $\widehat{f} := \tf_{x \rightarrow p} f$ in favour of readability.
		\item Let $A$ be the operator given by $A := i \, \partial_{xx}$ with domain $D(A) := H^3(\R) \subset H^1(R)$.
		\item For $V$ in $W^{1, \infty}(\R)$, let $B_V$ be the operator defined on $H^1(\R)$ by
		\begin{equation*}
			\forall \, f \in H^1(\R) \qquad (B_V f)(x) := -i \, V(x) \, f(x) \quad a.e \; .
		\end{equation*}
		In this case, the element $B_V f$ belongs to $H^1(\R)$; see for instance \citet[Chap. VI, Lem. 5.20]{en2000}.
	\end{enumerate}
\end{definition}

Let us prove now existence and uniqueness for the solution of equation \eqref{eq:schro-potential2} and provide a series representation.

\begin{theorem} \label{thm:existence_unique}
	Suppose that $u_0$ belongs to $H^3(\R)$ and that $V$ belongs to $W^{1, \infty}(\R)$. Then there exists a unique function $u : [0, +\infty) \longrightarrow H^3(\R) \subset H^1(\R)$ which is continuously differentiable with respect to the $H^1$-norm and which satisfies the Schrödinger equation \eqref{eq:schro-potential2}.\\
	Moreover the function $u$ can be represented as follows:
	\begin{equation*}
		\forall \, t \geqslant 0 \qquad \lim_{N \rightarrow + \infty} \left\| \, u(t) - \sum_{n = 1}^{N} S_n(t)u_0 \, \right\|_{H^1(\R)} = \, 0 \; ,
	\end{equation*}
	where
	\begin{equation*}
		\left\{ \begin{array}{rl}
			& \hspace{-2mm} S_1(t) u_0 := \tf_{p \rightarrow x}^{-1} \Big( e^{-i t p^2} \widehat{u_0}(p) \Big)  \\ [2mm]
			& \displaystyle \hspace{-2mm} S_{n+1}(t) u_0 := \int_0^t S_n(t-\s) \, B_V \, S_1(\s) u_0 \, d\s \quad , \quad \forall \, n \geqslant 1
		\end{array} \right. \; .
	\end{equation*}
\end{theorem}

\begin{remark} \label{rem:existence_unique}
	\begin{enumerate}
		\item The function $u: [0, +\infty) \longrightarrow H^1(\R)$ is called the \emph{classical solution} of the Schrödinger equation \eqref{eq:schro-potential2}; see \citet[Chap. II, Prop. 6.2]{en2000}. The series $\sum_{n \geqslant 1} S_n(t)u_0$ is called the \emph{Dyson-Phillips series} for the solution; we refer to \citet[Chap. III, Thm. 1.10]{en2000} for more details.
		\item For each $n \geqslant 1$, the term $S_{n+1}(t) u_0$ belongs at least to $H^1(\R)$ for all fixed $t \geqslant 0$ if $u_0 \in H^3(\R)$ and thus it defines a continuous function on $\R$. To evaluate it at any point $x \in \R$, let us define the evaluation operator $E_x$ on $H^1(\R)$ by
	\begin{equation*}
		\forall \, f \in H^1(\R) \qquad E_x f := f(x) \; .
	\end{equation*}
	Then $E_x$ is a bounded operator from $H^1(\R)$ into $\C$ thanks to the continuous embedding of $H^1(\R)$ into
	\begin{equation*}
		\mathcal{C}_0^0(\R):= \left\{ f \in \mathcal{C}^0(\R) \, \Big| \, \lim_{|x| \rightarrow +\infty} f(x) = 0 \right\} \; .
	\end{equation*}	
	Hence Proposition \ref{prop:bochner} is applicable and provides
	\begin{align}
		\forall \, x \in \R \qquad \big(S_{n+1}(t) u_0 \big)(x)	& = E_x \big(S_{n+1}(t) u_0 \big) \nonumber \\
			& = E_x \hspace{-1mm} \left( \int_0^t S_n(t-\s) \, B_V \, S_1(\s) u_0 \, d\s \right) \nonumber \\
			& \label{eq:evaluation} = \int_0^t \Big( S_n(t-\s) \, B_V \, S_1(\s) u_0 \Big)(x) \, d\s \; .
	\end{align}
	Note that the integral defining $S_{n+1}(t) u_0$ and the integral given in \eqref{eq:evaluation} are here interpreted as Bochner-integrals. In particular the integrand of \eqref{eq:evaluation} is complex-valued and, thanks to the construction of the Bochner-integral, it is actually an integral of Lebesgue-type. This evaluation process will be employed in Sections \ref{sec:s1} and \ref{sec:s2}.
	\end{enumerate}
\end{remark}

\begin{proof}[Proof of Theorem \ref{thm:existence_unique}]
	In order to apply results from semigroup theory, we start by rewri\-ting the Schrödinger equation \eqref{eq:schro-potential2} as an evolution equation of the form
	\begin{equation*} \label{eq:evol-eq}
		\left\{ \begin{array}{rl}
				& \hspace{-2mm} \dot{u}(t) = \big( A + B_V \big) u(t) \\ [2mm]
				& \hspace{-2mm} u(0) = u_0
		\end{array} \right. \; ,
	\end{equation*}
	where the operators $A$ and $B_V$ are given in Definition \ref{def:misc}.\\
	Now let us recall that the operator $\big(A, D(A) \big)$ is the generator of the strongly continuous semigroup $\big( T(t) \big)_{t \geqslant 0}$ on $H^1(\R)$ represented by
	\begin{equation} \label{eq:sgT}
		 \forall \, t \geqslant 0 \qquad T(t) f = \tf_{p \rightarrow x}^{-1} \Big( e^{-i t p^2} \widehat{f}(p) \Big) \; ,
	\end{equation}
	for $f \in H^1(\R)$. Moreover the operator $B_V$ belongs to $\mathcal{L}\big(H^1(\R) \big)$, the space of bounded operators from $H^1(\R)$ into itself; indeed we have
	\begin{align*}
		\big\| B_V f \big\|_{H^1(\R)}^2	& = \int_{\R} \big| V(x) f(x) \big|^2 \, dx +  \int_{\R} \big| (Vf)'(x) \big|^2 \, dx \\
											& = \int_{\R} \big| V(x) f(x) \big|^2 \, dx +  \int_{\R} \big| V'(x)f(x) \big|^2 \, dx + \int_{\R} \big| V(x)f'(x) \big|^2 \, dx \\
											& \leqslant 2 \, \big\| V \big\|_{W^{1, \infty}(\R)}^2 \big\| f \big\|_{H^1(\R)}^2 \; ,
	\end{align*}
	since $V \in W^{1, \infty}(\R)$. According to Proposition \ref{defprop:existence} and Theorem \ref{thm:perturb}, if $u_0 \in H^3(\R) = D(A)$ then the Schrödinger equation \eqref{eq:schro-potential2} has a unique classical solution belonging to $\mathcal{C}^1\big([0, +\infty), H^1(\R)\big)$. Moreover the solution is given by
	\begin{equation*}
		\forall \, t \geqslant 0 \qquad u(t) = S(t) u_0 \; ,
	\end{equation*}
	where $\big( S(t) \big)_{t \geqslant 0}$ is the semigroup generated by the operator $\big( A+B_V,D(A) \big)$ and belongs to $D(A) = H^3(\R)$ for all $t \geqslant 0$.\\
	Employing now Theorem \ref{thm:dyson}, the solution of equation \eqref{eq:schro-potential2} can be represented as follows,
	\begin{equation*}
		\forall \, t \geqslant 0 \qquad \lim_{N \rightarrow + \infty} \left\| \, u(t) - \sum_{n = 1}^{N} S_n(t)u_0 \, \right\|_{H^1(\R)} = \, 0 \; ,
	\end{equation*}
	where $S_1(t) := T(t)$ and
	\begin{equation*}
		S_{n+1}(t) u_0 := \int_0^t S_n(t-\s) \, B_V \, T(\s) u_0 \, ds \; .
	\end{equation*}
	According to equality \eqref{eq:sgT}, we have $S_1(t)u_0 = T(t) u_0 = \tf_{p \rightarrow x}^{-1} \big( e^{-i t p^2} \widehat{u_0}(p) \big)$, which ends the proof.
\end{proof}

\section{Potentials in bounded Fourier-frequency bands} \label{sec:potential}

The main goal of this short section is to introduce and to illustrate the hypotheses on the potential. Roughly speaking, we assume the potential to be in a bounded frequency band in order to isolate the effect of the frequencies of the potential on the time-asymptotic motion of the Dyson-Phillips series terms. For the sake of illustration, we provide also an explicit family of potentials which verify this frequency band hypothesis.\\

The hypotheses of interest are given in the following condition.\\

\noindent \textbf{Condition} ($\mathcal{P}_{[a,b]}^k$). Let $k \in \N$ and let $a < b$ be two finite real numbers. \\
An element $V$ of $L^2(\R)$ satisfies Condition ($\mathcal{P}_{[a,b]}^k$) if and only if $\widehat{V}$ is a $\mathcal{C}^k$-function on $\R$ which verifies $\supp \, \widehat{V} \subseteq [a,b]$.
	
\begin{remark} \label{rem:potential_hyp}
	\begin{enumerate}
		\item The set of functions satisfying Condition ($\mathcal{P}_{[a,b]}^k$) is non-empty. Indeed if $U$ is a $\mathcal{C}^k$-function supported on $[a,b]$ then it belongs to $L^2(\R)$. Hence, thanks to the fact that $\tf_{x \rightarrow p}$ is a bijective map from $L^2(\R)$ onto itself, there exists $V$ in $ L^2(\R)$ such that $U = \widehat{V}$. In particular $V$ satisfies Condition ($\mathcal{P}_{[a,b]}^k$).
		\item Under Condition ($\mathcal{P}_{[a,b]}^k$), a potential $V$ is actually analytic because of the boundedness of the support of $\widehat{V}$.
		\item If a potential $V$ verifies Condition ($\mathcal{P}_{[a,b]}^k$), with $k \geqslant 1$, then it is bounded on $\R$ as well as its first derivative. Indeed the functions $p \longmapsto \widehat{V}(p)$ and $p \longmapsto p \, \widehat{V}(p)$ belong to $L^1(\R)$ under Condition ($\mathcal{P}_{[a,b]}^k$), which permits to bound the norms $\| V \|_{L^\infty(\R)}$ and $\| V' \|_{L^\infty(\R)}$. In particular, the potential $V$ belongs to $W^{1,\infty}(\R)$ and so the associated operator $B_V$ defined in Section \ref{sec:dyson} belongs to $\mathcal{L}\big( H^1(\R) \big)$.
		\item A potential $V$ satisfying Condition ($\mathcal{P}_{[a,b]}^k$) is not necessarily real-valued: it is real-valued if and only if $\widehat{V}$ verifies
	\begin{equation*}
		\forall \, p \in \R \qquad \widehat{V}(-p) = \overline{\widehat{V}(p)} \; .
\end{equation*}			
		In this case, the support of $\widehat{V}$ is contained in a symmetric interval centred on the origin; hence we have especially $a = -b$.
	\end{enumerate}
\end{remark}

To illustrate the above Condition ($\mathcal{P}_{[a,b]}^k$), we construct a family of admissible potentials. It is noteworthy that any element of this family is approximately localised in space around a given point with arbitrary precision. This is a direct consequence of a more generic result stated in Lemma \ref{lem:localisation_fourier}; see Section \ref{sec:appendix}.

\begin{proposition} \label{prop:admissible_potentials}
	Let $k \geqslant 1$ be an integer, let $a$, $b$ and $x_0$ be three finite real numbers such that $a < b$, and let $v$ be a $\mathcal{C}^k$-function such that $\supp \, v \subseteq [-1,1]$.
	Let $V$ be the element of $L^2(\R)$ whose Fourier transform $\widehat{V}$ is the complex-valued function given by
	\begin{equation*}
		\forall \, p \in \R \qquad \widehat{V}(p) := v \hspace{-1mm} \left( \frac{2p - (a+b)}{b-a} \right) e^{-i x_0 p} \; .
	\end{equation*}
	Then $V$ verifies Condition ($\mathcal{P}_{[a,b]}^k$) and satisfies
	\begin{equation*}
		\forall \, c > 0 \qquad \int_{|x-x_0| \geqslant c} \big| V(x) \big|^2 dx \leqslant \frac{2}{c^2} \, \frac{1}{b-a} \, \big\| v' \big\|_{L^2(\R)}^2 \; .
	\end{equation*}
\end{proposition}

\begin{proof}
	Straightforward application of Lemma \ref{lem:localisation_fourier} to the function $V$.
\end{proof}

\section{Time-asymptotic behaviour of the free solution} \label{sec:s1}

This section is devoted to the study of the time-asymptotic propagation of the first term of the Dyson-Phillips series, which is actually the free Schrödinger solution according to Theorem \ref{thm:existence_unique}. More precisely we show that the first term tends to be time-asymptotically localised in a space-time cone if the initial datum is in a frequency band. To do so, we apply the method developed in \citet{amhdr2012, amd2017} and refined in \citet{dewez2020}, whose main result is slightly extended in Theorem \ref{thm:oi_exp} below. This result will be also exploited to study the second term of the Dyson-Phillips series but extra technical arguments will be required; we refer to Section \ref{sec:s2} for the study of the second term.\\

We start by recalling the definition of a space-time cone related to the (frequency) interval $[\widetilde{p}_1, \widetilde{p}_2]$ with origin $(t_0, x_0) \in \R \times \R$.

\begin{definition}
	Let $t_0$, $x_0$, $\tilde{p}_1$ and $\tilde{p}_2$ be four finite real numbers with $\tilde{p}_1 < \tilde{p}_2$.
	\begin{enumerate}
		\item We define the space-time cone $\mathfrak{C}\big( [\tilde{p}_1, \tilde{p}_2], (t_0, x_0) \big)$ as follows:
		\begin{equation} \label{eq:cone}
			\mathfrak{C}\big( [\tilde{p}_1, \tilde{p}_2], (t_0, x_0) \big) := \left\{ (t,x) \in \big( [0,+\infty) \backslash \{t_0 \} \big) \times \R \, \bigg| \, 2 \, \tilde{p}_1 \leqslant \frac{x - x_0}{t - t_0} \leqslant 2 \, \tilde{p}_2 \right\} \; .
		\end{equation}
		\item Let $\mathfrak{C}\big( [\tilde{p}_1, \tilde{p}_2], (t_0, x_0) \big)^c$ be the complement of the space-time cone $\mathfrak{C}\big( [\tilde{p}_1, \tilde{p}_2], (t_0, x_0) \big)$ in $\big( [0,+\infty) \backslash \{t_0 \} \big) \times \R$ .
	\end{enumerate}	 
\end{definition}

\begin{remark}
	When one considers general dispersive equations of the form \eqref{eq:disp_equation}, the above cone is defined by the following inequalities:
	\begin{equation*}
		f'(\tilde{p}_1) \leqslant \frac{x - x_0}{t - t_0} \leqslant f'(\tilde{p}_2) \; .
	\end{equation*}
	This explains the factors $2$ appearing in \eqref{eq:cone} since the symbol in the Schrödinger setting is $f(p) = p^2$.
\end{remark}

In the rest of this paper, the initial datum will be assumed to satisfy the following frequency band hypothesis.\\

\noindent \textbf{Condition} ($\mathcal{I}_{[p_1, p_2]}^\ell$). Let $\ell \in \N$ and let $p_1 < p_2$ be two finite real numbers. \\
	An element $u_0$ of $H^3(\R)$ satisfies Condition ($\mathcal{I}_{[p_1, p_2]}^\ell$) if and only if $\widehat{u_0}$ is a $\mathcal{C}^\ell$-function on $\R$ which verifies $\supp \, \widehat{u_0} \subseteq [p_1,p_2]$.
	
\begin{remark} \label{rem:i_p1p2}
	\begin{enumerate}
		\item Following arguments similar to those of Remark \ref{rem:potential_hyp}, one proves that the set of elements of $H^3(\R)$ satisfying Condition ($\mathcal{I}_{[p_1, p_2]}^\ell$) is non-empty and that it is a set of analytic functions.
		\item If $u_0$ satisfies Condition ($\mathcal{I}_{[p_1, p_2]}^\ell$) then its Fourier transform is an integrable function on $\R$. Hence for $t \geqslant 0$, the first term $S_1(t) u_0$ of the Dyson-Phillips series introduced in Theorem \ref{thm:existence_unique} is actually a complex-valued function on $\R$ given by
		\begin{equation} \label{eq:s1}
			\forall \, x \in \R \qquad \big( S_1(t) u_0 \big)(x) = \frac{1}{2 \pi} \int_{p_1}^{p_2} \widehat{u_0}(p) \, e^{-itp^2 + ixp} \, dp \; .
		\end{equation}
	\end{enumerate}
\end{remark}

Let us now talk briefly about \citet[Thm. 1.1]{dewez2020} and motivate the need for a slight extension. The method used to obtain \citet[Thm. 1.1]{dewez2020} is entirely based on the integral representation \eqref{eq:s1} of the free Schrödinger solution. The proof consists in firstly making a space-time shift in the above wave packet \eqref{eq:s1}, secondly in factorising the phase function by the time to get a generic oscillatory integral and finally in applying carefully a stationary phase method; see \citet[Thm. 3.3 and 3.4]{dewez2020}. We mention that these two results are adaptations of the version of the stationary phase method given in \cite{amd2017}, which is itself a modern and refined version of the results established in \citet[Sec. 2.9]{erdelyi1956}. This adapted version in \citet{dewez2020} provides a remainder estimate which is explicit with respect to the space-time shift parameter; note that this parameter describes actually the origin of the cone. This flexibility allows then to choose the shift parameter which both minimises the remainder estimate and makes stable propagation features under time-asymptotic expansions.

However we need here to extend slightly \citet[Thm. 1.1]{dewez2020} because it has been originally proved for free Schrödinger wave packets with initial data in the Schwartz space, which is not necessarily the case in the present paper. We mention that the choice for the Schwartz space in \citet{dewez2020} has been done for the sake of clarity but a careful look at the proof shows that the smoothness and decay assumptions on the initial datum can be relaxed. Further, in view of the application to two the second of the Dyson-Phillips series, we extend also \citet[Thm. 1.1]{dewez2020} to wave packets of the form 
\begin{equation*}
	\forall \, (t,x) \in \R \times \R \qquad J_U(t,x) := \frac{1}{2 \pi} \int_\R U(p) \, e^{-itp^2 + ixp} \, dp \; ,
\end{equation*}
where the compactly supported amplitude $U$ is not required to be the Fourier transform of an initial datum $u_0$.

\begin{theorem} \label{thm:oi_exp}
	Let $\ell \geqslant 1$ be an integer, let $p_1$, $p_2$, $\widetilde{p}_1$ and $\widetilde{p}_2$ be four finite real numbers such that $[p_1, p_2] \subset \big(\widetilde{p}_1, \widetilde{p}_2 \big)$. Let $U$ be a $\mathcal{C}^\ell$-function such that 
	\begin{equation*}
		\supp \, U \subseteq [p_1, p_2] \qquad \text{and} \qquad \frac{1}{\sqrt{2 \pi}} \, \| U \|_{L^2(\R)} = 1 \; .
\end{equation*}		
	And we define
	\begin{align*}
		& \bullet \quad t^* = \argmin_{\tau \in \R} \left( \int_\R x^2 \, \big| J_U(t,x) \big|^2 \, dx - \Big( \int_\R x \, \big| J_U(t,x) \big|^2 \, dx \Big)^2 \right) ; \\
		& \bullet \quad x^* = \int_\R x \, \big| J_U(t^*,x) \big|^2 \, dx \; .
	\end{align*}
	Then for all $\displaystyle (t,x) \in \mathfrak{C}\big( [\widetilde{p}_1, \widetilde{p}_2], (t^*, x^*) \big)$, we have
	\begin{align*}
		& \left| J_U(t,x) - \frac{1}{\sqrt{2 \pi}} \, e^{- sgn(t-t^*) i \frac{\pi}{4}} \, e^{-it \big(\frac{x-x^*}{t - t^*}\big)^2 + ix \frac{x-x^*}{t - t^*}} \, U \left( \frac{x-x^*}{t - t^*} \right) |t-t^*|^{-\frac{1}{2}} \right| \nonumber \\
		 & \hspace{1cm} \leqslant C_1(\delta, \tilde{p}_1, \tilde{p}_2) \, \sqrt{\int_\R x^2 \, \big| J_U(t^*,x) \big|^2 \, dx - \Big(\int_\R x \, \big| J_U(t^*,x) \big|^2 \, dx \Big)^2} \, |t-t^*|^{-\delta} \; ,
	\end{align*}
	where the real number $\delta$ is arbitrarily chosen in $\big( \frac{1}{2}, \frac{3}{4} \big)$. And for all $\displaystyle (t,x) \in \mathfrak{C}\big( [\widetilde{p}_1, \widetilde{p}_2], (t^*, x^*) \big)^c$, we have
	\begin{align*}
		\Big| J_U(t,x) \Big|
			& \leqslant \Bigg( C_2(p_1, p_2, \tilde{p}_1, \tilde{p}_2) \, \sqrt{\int_\R x^2 \, \big| J_U(t^*,x) \big|^2 \, dx - \Big( \int_\R x \, \big| J_U(t^*,x) \big|^2 \, dx \Big)^2} \\
			& \hspace{1.5cm} +  C_3(p_1, p_2, \tilde{p}_1, \tilde{p}_2) \, \big\| U \big\|_{L^\infty(\R)} \Bigg) \, |t-t^*|^{-1} \; .
	\end{align*}
	All the above constants are explicitly given in \citet[Thm. 1.1]{dewez2020}.
\end{theorem}

\begin{proof}
	In \citet{dewez2020}, Theorem 1.1 is a straightforward consequence of Corollary 2.5, which is itself an application of the more general result Theorem 2.2. Following the proof of this theorem, we rewrite the wave packet $J_U(t,x)$ as follows:
	\begin{equation*}
		J_U(t,x) = \int_\R \frac{1}{2 \pi} \, U(p) \, e^{-i t_0 p^2 + i x_0 p} \, e^{i (t-t_0) \big(\frac{x-x_0}{t-t_0} \, p - p^2 \big)} \, dp \; ,
	\end{equation*}		
	where $t_0, x_0 \in \R$, and we apply Theorems 3.3 and 3.4 from \citet{dewez2020} to the above rewriting. These two theorems provide asymptotic expansions and uniform estimates of generic oscillatory integrals. They require only the amplitude of the integral to be con\-ti\-nuous\-ly differentiable on $\R$ with compact support, which is the case in the present setting. This provides then a time-asymptotic expansion of $J_U(t,x)$ and a uniform estimate respectively inside and outside the cone $\mathfrak{C}\big( [\widetilde{p}_1, \widetilde{p}_2], (t_0, x_0) \big)$. As in the end of the proof of Corollary 2.5 from \citet{dewez2020}, we finish by setting the origin of the cone: $t_0 = t^*$, $x_0 = x^*$, where $t^*, x^* \in \R$ are defined in the statement of the theorem.\\
	 It is noteworthy that the existence of $t^*$ and $x^*$ depends on the fact that the (spatial) variance of $\big| J_U(t_0,.) \big|^2$ is well-defined for any time $t_0 \in \R$. To prove the existence of the variance, we remark that the present assumptions on the amplitude $U$ are actually sufficient so that $\big| J_U(t_0,.) \big|^2$ has a moment of order 2 about any $x_0 \in \R$ for all $t_0 \in \R$. Indeed Plancherel theorem and standard properties of the Fourier transform lead to
	 \begin{equation*}
	 	\Big\| x \longmapsto (x-x_0) \, J_U(t_0,x_0) \Big\|_{L^2(\R)}^2 = \frac{1}{2 \pi} \, \Big\| \partial_p\Big(p \longmapsto U(p) \, e^{-i t_0 p^2 + i x_0 p} \Big) \Big\|_{L^2(\R)}^2 \; ,
	 \end{equation*}
	and the right-hand side is well-defined for all $t_0, x_0 \in \R$. In particular, the variance of $\big| J_U(t_0,.) \big|^2$ is well-defined for any time $t_0 \in \R$.\\
	Finally we emphasise that the $L^\infty$-norm of $U$ appearing in the second inequality of the theorem is not bounded by the $L^1$-norm of $\tf_{p \mapsto x}^{-1} U$ as it is done in \citet{dewez2020}. Indeed we are not in position to claim that the $L^2$-element $\tf_{p \mapsto x}^{-1} U$ is actually an integrable function so we do not apply here the classical inequality $\| U \|_{L^\infty(\R)} \leqslant \| \tf_{p \mapsto x}^{-1} U \|_{L^1(\R)}$.
\end{proof}

For initial data satisfying Condition ($\mathcal{I}_{[p_1, p_2]}^\ell$), a straightforward application of the preceding theorem leads to a time-asymptotic expansion of the term $S_1(t) u_0$ inside the cone $\mathfrak{C}\big( [\widetilde{p}_1, \widetilde{p}_2], (t_1^*, x_1^*) \big)$ and a uniform estimate outside, where $[p_1, p_2] \subset \big(\widetilde{p}_1, \widetilde{p}_2 \big)$ and $t_1^*$ and $x_1^*$ are defined below. This result shows that $S_1(t) u_0$ is time-asymptotically close to a term which is supported on the time-dependent interval $\big[x_1^* + 2 \, \widetilde{p}_1 \,(t-t_1^*), x_1^* + 2 \, \widetilde{p}_2 \,(t-t_1^*)\big]$, providing then information on the motion and the dispersion of $S_1(t) u_0$ for $t$ far from $t_1^*$.

\begin{corollary} \label{cor:s1_exp}
	Let $\ell \geqslant 1$ be an integer, let $p_1$, $p_2$, $\widetilde{p}_1$ and $\widetilde{p}_2$ be four finite real numbers such that $[p_1, p_2] \subset \big(\widetilde{p}_1, \widetilde{p}_2 \big)$. Suppose that $u_0$ satisfies Condition ($\mathcal{I}_{[p_1, p_2]}^\ell$) with $\| u_0 \|_{L^2(\R)} = 1$ and define
	\begin{align*}
		& \bullet \quad t_1^* = \argmin_{\tau \in \R} \left( \int_\R x^2 \, \big| S_1(\s)u_0(x) \big|^2 \, dx - \Big( \int_\R x \, \big| S_1(\s)u_0(x) \big|^2 \, dx \Big)^2 \right) ; \\
		& \bullet \quad x_1^* = \int_\R x \, \big| S_1(t_1^*)u_0(x) \big|^2 \, dx \; .
	\end{align*}
	Then for all $\displaystyle (t,x) \in \mathfrak{C}\big( [\widetilde{p}_1, \widetilde{p}_2], (t_1^*, x_1^*) \big)$, we have
	\begin{align*}
		& \left| S_1(t)u_0(x) - \frac{1}{\sqrt{2 \pi}} \, e^{- sgn(t-t_1^*) i \frac{\pi}{4}} \, e^{-it \big(\frac{x-x_1^*}{t - t_1^*}\big)^2 + ix \frac{x-x_1^*}{t - t_1^*}} \, \widehat{u}_0 \left( \frac{x-x_1^*}{t - t_1^*} \right) |t-t_1^*|^{-\frac{1}{2}} \right| \nonumber \\
		 & \hspace{1cm} \leqslant C_1(\delta, \tilde{p}_1, \tilde{p}_2) \, \sqrt{\int_\R x^2 \, \big| S_1(t_1^*)u_0(x) \big|^2 \, dx - \Big(\int_\R x \, \big| S_1(t_1^*)u_0(x) \big|^2 \, dx \Big)^2} \, |t-t_1^*|^{-\delta} \; , \label{eq:schro_ae}
	\end{align*}
	where the real number $\delta$ is arbitrarily chosen in $\big( \frac{1}{2}, \frac{3}{4} \big)$. And for all $\displaystyle (t,x) \in \mathfrak{C}\big( [\widetilde{p}_1, \widetilde{p}_2], (t_1^*, x_1^*) \big)^c$, we have
	\begin{align*}
		\Big| S_1(t)u_0(x) \Big|
			& \leqslant \Bigg( C_2(p_1, p_2, \tilde{p}_1, \tilde{p}_2) \, \sqrt{\int_\R x^2 \, \big| S_1(t_1^*)u_0(x) \big|^2 \, dx - \Big( \int_\R x \, \big| S_1(t_1^*)u_0(x) \big|^2 \, dx \Big)^2} \\
			& \hspace{1.5cm} +  C_3(p_1, p_2, \tilde{p}_1, \tilde{p}_2) \, \big\| \widehat{u}_0 \big\|_{L^\infty(\R)} \Bigg) \, |t-t_1^*|^{-1} \; .
	\end{align*}
	All the above constants are explicitly given in \citet[Thm. 1.1]{dewez2020}.
\end{corollary}

\begin{proof}
	It is sufficient to apply Theorem \ref{thm:oi_exp} to the oscillatory integral \eqref{eq:s1}, which represents the term $S_1(t) u_0$, with $U = \widehat{u}_0$. Condition ($\mathcal{I}_{[p_1, p_2]}^\ell$) assures in particular that the hypotheses of Theorem \ref{thm:oi_exp} are verified.
\end{proof}

\section{Time-asymptotic behaviour of the wave packet issued by a first interaction with the potential} \label{sec:s2}

In this section, we focus on the time-asymptotic propagation of the second term $S_2(t) u_0$ of the Dyson-Phillips series. As in the preceding section, the aim is to expand this term in a well-suited cone to reflect its time-asymptotic localisation in space. In particular, the frequency band hypothesis on the potential will highlight the influence of the frequencies of the potential on the propagation of $S_2(t) u_0$.

From a technical point of view, we show in the first subsection that the second term $S_2(t) u_0$ can be written as a wave packet with time-dependent amplitude. We prove then in the second subsection that $S_2(t) u_0$ tends to be time-asymptotically close to a free wave packet, \emph{i.e.} a wave packet with an amplitude independent from time. To do so, we show that the time-dependent amplitude converges to a limit amplitude. This step is required because if we applied Theorem \ref{thm:oi_exp} directly to $S_2(t) u_0$, then the (absolute value of the) coefficient of the resulting first term of the expansion would depend on time, preventing from deriving the time-decay rate. In the last subsection, we deduce the time-asymptotic behaviour of $S_2 (t) u_0$ by applying Theorem \ref{thm:oi_exp} to the free wave packet with the limit amplitude, providing finally a time-expansion of $S_2 (t) u_0$ in a cone and a uniform estimate outside.

As explained in the preceding subsection and in the introduction, the cone inclination gives information on both the motion and the dispersion. Here it is noteworthy that the inclination of the cone associated with $S_2(t) u_0$ depends explicitly on the frequencies of both the initial datum and the potential, highlighting especially the influence of the frequencies of the potential on the propagation of this perturbed term.

\subsection{Wave packet representation for the second term of the Dyson-Phillips series}

In the only proposition of this subsection, we give a representation of the second term $S_2(t) u_0$ as an wave packet with a time-dependent amplitude.\\
The approach to prove Proposition \ref{prop:rewriting} is based on the explicit formula of $S_2(t) u_0$ given in Theorem \ref{thm:existence_unique} and on applications of basic properties of the Fourier transform as well as Fubini's theorem. Let us note in particular that the amplitude of the resulting wave packet is actually supported on the same bounded interval for any time.

\begin{proposition} \label{prop:rewriting}
	Let $k \geqslant 1$ and $\ell \geqslant 0$ be two integers and suppose that $u_0 \in H^3(\R)$ and $V \in L^2(\R)$ satisfy respectively Conditions \emph{($\mathcal{I}_{[p_1,p_2]}^\ell$)} and \emph{($\mathcal{P}_{[a,b]}^k$)}. Let $t \geqslant 0$ and let $W(t,.) : \R \longrightarrow \C$ be the function defined by
	\begin{equation*}
		W(t,p) := -i \int_0^t \widetilde{W}(\s,p) \, e^{i \s p^2} \, d\s \; ,
	\end{equation*}
	where we have defined the function $\widetilde{W}: [0,t] \times \R \longrightarrow \C$ as follows:
	\begin{align*}
		\widetilde{W}(\s,p) := \Big( \widehat{V} \ast \big( e^{-i \s \, \cdot^2} \, \widehat{u_0}(.) \big) \Big)(p) = \int_{a}^{b} \widehat{V}(y) \, \widehat{u_0}(p-y) \, e^{- i \s (p-y)^2} \, dy \; .
	\end{align*}		
	Then the support of $W(t,.)$ is contained in $[p_1 + a, p_2 + b]$ for all $t \geqslant 0$, and
	\begin{equation*}
		\forall \, (t,x) \in [0, +\infty) \times \R \qquad \big( S_2(t) u_0 \big)(x) = \frac{1}{2 \pi} \int_\R W(t,p) \, e^{- i t p^2 + i x p} \, dp \; .
	\end{equation*}
\end{proposition}

\begin{proof}
	Let $t \geqslant 0$ and $x \in \R$. By using equality \eqref{eq:evaluation} from Remark \ref{rem:existence_unique} and the Fourier representation of $S_1(\s)u_0 = T(\s)u_0$, we obtain
	\begin{align}
		\big(S_2(t) u_0 \big)(x)	& = \int_0^t \big( T(t-\s) \, B \, T(\s) u_0 \big)(x) \, d\s \nonumber \\
					& = -i \int_0^t \tf_{p \rightarrow x}^{-1} \bigg( e^{-i (t-\s) \, p^2} \tf_{x \rightarrow p} \Big( V(x) \tf_{p \rightarrow x}^{-1} \big( e^{-i \s \, p^2} \, \widehat{u_0}(p) \big)(x) \Big)(p) \bigg)(x) \, d\s \nonumber \\
					& = -i \int_0^t \tf_{p \rightarrow x}^{-1} \bigg( e^{-i (t-\s) \, p^2} \Big( \widehat{V} \ast \big( e^{-i \s \, \cdot^2} \, \widehat{u_0}(.) \big) \Big)(p) \bigg)(x) \, d\s \nonumber \\
					& = -i \int_0^t \tf_{p \rightarrow x}^{-1} \Big( e^{-i (t-\s) \, p^2} \, \widetilde{W}(\s, \, p) \Big)(x) \, d\s  \label{eq:tildeW} \; .
	\end{align}
	Now let us remark that, for any $\s \in [0,t]$, the function $\widetilde{W}(\s,.)$ is the convolution of two compactly supported and continuous functions on $\R$. Therefore $\widetilde{W}(\s,.)$ is also a continuous function on $\R$ such that
	\begin{equation*}
		\supp \, \widetilde{W}(\s,.) = \supp \, \Big( \widehat{V} \ast \big( e^{-i \s \, \cdot^2} \, \widehat{u_0}(.) \big) \Big) \subseteq \overline{[a,b] + [p_1 , p_2]} = [p_1 + a , p_2 + b] \; ,
	\end{equation*}
	for any $\s \in [0,t]$. In particular, we deduce that the support of $W(t,.)$ is also contained in $[p_1 + a , p_2 + b]$. It follows that $\widetilde{W}(\s,.)$ is an integrable function and so the quantity $\tf_{p \rightarrow x}^{-1} \big( e^{-i (t-\s) \, p^2} \, \widetilde{W}(\s, \, p) \big)(x)$ can be given by the integral representation of the inverse Fourier transform for integrable functions with respect to the variable $p$, \emph{i.e.}
	\begin{equation*}
		\tf_{p \rightarrow x}^{-1} \Big( e^{-i (t-\s) \, p^2} \, \widetilde{W}(\s, p) \Big)(x) = \frac{1}{2 \pi} \int_\R \widetilde{W}(\s,p) \, e^{-i (t-\s) p^2 + i x p} \, dp \; .
	\end{equation*}
	Combining this with equality \eqref{eq:tildeW} leads to
	\begin{equation*}
		\big( S_2(t) u_0 \big)(x) = -i \int_0^t \left(\frac{1}{2 \pi} \int_\R \widetilde{W}(\s,p) \, e^{-i (t-\s) p^2 + i x p} \, dp \right) d\s \; .
	\end{equation*}
	Since the integrand in the preceding double integral is a continuous function on the compact domain $[0,t] \times [p_1 + a, p_2 + b]$, we can apply Fubini's theorem to obtain the desired equality, namely,
	\begin{align*}
		\big( S_2(t) u_0 \big)(x)	& = \frac{1}{2 \pi} \int_\R \left(-i \int_0^t \widetilde{W}(\s,p) \, e^{i \s p^2} \, d\s \right) e^{-i t p^2 + i x p} \, dp \\
										& = \frac{1}{2 \pi} \int_\R W(t,p) \, e^{-i t p^2 + i x p} \, dp \; .
	\end{align*}
\end{proof}

\begin{remark}
	In the preceding result, the Fourier transform $\widehat{V}$ is assumed to be at least a $\mathcal{C}^1$-function to assure that the operator $B_V$ introduced in Definition \ref{def:misc} belongs to $\mathcal{L}\big( H^1(\R) \big)$ as explained in Remark \ref{rem:potential_hyp} 3.
\end{remark}


\subsection{Limit of the time-dependent amplitude}

In this subsection, we are interested in the limit as $t$ tends to infinity of the time-dependent amplitude $W(t,.)$ defined in Proposition \ref{prop:rewriting}. This limit is here proved to be a $\mathcal{C}^1$-function with compact support and an upper bound for the convergence speed is provided.\\

We begin with the study of the sum of two parametric integrals which will be proved to be the limit of $W(t,.)$ in Proposition \ref{prop:exp_winfty}. In the following lemma, we aim at showing that this sum defines a continuously differentiable function with support contained in the interval $[p_1 + a, p_2 + b]$. We also provide an explicit estimate of one of the two parametric integrals which will be used to derive an upper bound for the convergence speed.\\
From a technical point of view, we assume a certain regularity for the Fourier transforms of the initial datum $u_0$ and of the potential $V$. This additional assumption is exploited together with the frequency band hypotheses to apply classical results on parametric integrals. Let us also emphasise that we assume the frequency $0$ to be outside the support of $\widehat{u_0}$. This assumption seems to be a key point in our approach since it assures that the integrands defining the parametric integrals do not have singular points.

\begin{lemma} \label{lem:prop_winfty}
	Let $k \geqslant 3$ and $\ell \geqslant 4$ be two integers and suppose that $u_0 \in H^3(\R)$ and $V \in L^2(\R)$ satisfy respectively Conditions \emph{($\mathcal{I}_{[p_1,p_2]}^\ell$)} and \emph{($\mathcal{P}_{[a,b]}^k$)}. Assume in addition that $0 \notin [p_1, p_2]$. Let $s > 0$, $p \in \R$ and define
	\begin{align*}
		& \bullet \quad \displaystyle W_{\infty, 1}^s(p) := -i \int_0^s \int_a^b \widehat{V}(y) \, \widehat{u_0}(p-y) \, e^{-i \s (p-y)^2} \, dy \, e^{i \s p^2} d\s \; ; \\
		& \bullet \quad \displaystyle W_{\infty, 2}^s(p) := \frac{1}{8} \int_s^\infty \int_a^b \partial_y \left[ \frac{1}{p-y} \partial_y \left[ \frac{1}{p-y} \partial_y \frac{\widehat{V}(y) \, \widehat{u_0}(p-y)}{p-y} \right] \right] e^{-i \s (p-y)^2} dy \, \s^{-3} e^{i \s p^2} d\s \; ; \\
		& \bullet \quad \displaystyle W_\infty^s(p) := W_{\infty, 1}^s(p) + W_{\infty, 2}^s(p) \; .
	\end{align*}
	Then
	\begin{enumerate}
		\item the sum $W_\infty^s$ defines a continuously differentiable function on $\R$ with support contained in $[p_1 + a, p_2 + b]$;
		\item for all $p \in \R$,
		\begin{equation} \label{eq:w2_ineq}
			\big| W_{\infty, 2}^s(p) \big| \leqslant r(p) \, s^{-2} \; ,
		\end{equation}
		where $r: \R \longrightarrow \R$ is an integrable function defined in \eqref{eq:def_r}.
	\end{enumerate}
\end{lemma}

\begin{proof}
	Let $s > 0$ and $p \in \R$.
	\begin{enumerate}
		\item From the definition of $W(t,p)$ given in Proposition \ref{prop:rewriting}, where $t \geqslant 0$, we observe that $W_{\infty, 1}^s (p) = W(s,p)$. From this, we deduce that $W_{\infty, 1}^s$ defines a function on $\R$ with support contained in $[p_1 + a, p_2 + b]$. Furthermore $W_{\infty, 1}^s(p)$ is a parameter-dependent $(\s, y)$-integral of a $\mathcal{C}^1$-function (actually a $\mathcal{C}^4$-function) with respect to $p$ with integration domain given by $[0, s] \times [a, b]$. So it defines a $\mathcal{C}^1$-function on $\R$.
		
		Let us now show that $W_{\infty, 2}^s$ defines also a $\mathcal{C}^1$-function on $\R$ whose support is contained in $[p_1 + a, p_2 + b]$. Let us first prove that it is well-defined for all $p \in \R$ by showing that its integrand
		\begin{equation*}
			f(p, \s, y) := \frac{1}{8} \, \partial_y \left[ \frac{1}{p-y} \partial_y \left[ \frac{1}{p-y} \partial_y \frac{\widehat{V}(y) \, \widehat{u_0}(p-y)}{p-y} \right] \right] e^{-i \s (p-y)^2} \s^{-3} e^{i \s p^2}
		\end{equation*}
		is absolutely integrable with respect to $(\s, y)$. To this end, we observe that $f(p, \s, y) = 0$ on the line $y = p$ since $\widehat{u}_0(0) = 0$ by hypothesis on the support of $\widehat{u}_0$; we deduce that $f(p, \s, y)$ is well-defined for all $p \in \R$, $y \in [a,b]$ and $\s \geqslant s$. Further we have
		\begin{equation*}
			\big| f(p, \s, y) \big| \leqslant\left| \frac{1}{8} \, \partial_y \left[ \frac{1}{p-y} \partial_y \left[ \frac{1}{p-y} \partial_y \frac{\widehat{V}(y) \, \widehat{u_0}(p-y)}{p-y} \right] \right] \right|\s^{-3} \; ;
		\end{equation*}
		note that the right-hand side is a $(\s, y)$-integrable function on $[s, +\infty) \times [a,b]$ since in particular the functions $\widehat{V}$ and $\widehat{u}_0$ as well as their three first derivatives are continuous. Hence the function $f(p, ., .)$ is absolutely integrable and $W_{\infty, 2}^s(p)$ is well-defined.\\
		Let us now study the support of $W_{\infty, 2}^s: \R \longrightarrow \C$. We note that for fixed $\s \geqslant s$, $f(p, \s, .)$ is equal to zero outside the interval
		\begin{equation*}
			I_p := \big\{ y \in [a,b] \, \big| \, p-y \in [p_1, p_2] \big\} = [a,b] \cap [p - p_2 , p - p_1] \; ,
		\end{equation*}
		since $\widehat{V}$ and $\widehat{u}_0$ are respectively supported on $[a,b]$ and $[p_1, p_2]$. The interval $I_p$ is empty for $p \notin [p_1 + a, p_2 + b]$ so, in this case, the integral defining $W_{\infty, 2}^s(p)$ is equal to 0 proving that
		\begin{equation*}
			\supp \, W_{\infty, 2}^s \subseteq [p_1 + a, p_2 + b] \; .
		\end{equation*}
		We finish by proving that $W_{\infty, 2}^s$ is a $\mathcal{C}^1$-function on $\R$. The hypothesis on the regularity of $\widehat{u}_0$ permits especially to show that $f(., \s, y)$ is a $\mathcal{C}^1$-function for all $y \in [a,b]$ and $\s \geqslant s$. Further we have
		\begin{align*}
			\big| \partial_p f(p, \s, y) \big|
				& \leqslant \left| \frac{1}{8} \, \partial_p \partial_y \left[ \frac{1}{p-y} \partial_y \left[ \frac{1}{p-y} \partial_y \frac{\widehat{V}(y) \, \widehat{u_0}(p-y)}{p-y} \right] \right] \right|\s^{-3} \\
				& \hspace{5mm} + \left| \frac{1}{4} \, y \, \partial_y \left[ \frac{1}{p-y} \partial_y \left[ \frac{1}{p-y} \partial_y \frac{\widehat{V}(y) \, \widehat{u_0}(p-y)}{p-y} \right] \right] \right|\s^{-2} \\
				& \leqslant \left\| (p, y) \mapsto \frac{1}{8} \, \partial_p \partial_y \left[ \frac{1}{p-y} \partial_y \left[ \frac{1}{p-y} \partial_y \frac{\widehat{V}(y) \, \widehat{u_0}(p-y)}{p-y} \right] \right] \right\|_{L^\infty(\R \times [a,b])}\s^{-3} \\
				& \hspace{5mm} + \left\| (p, y) \mapsto \frac{1}{4} \, y \, \partial_y \left[ \frac{1}{p-y} \partial_y \left[ \frac{1}{p-y} \partial_y \frac{\widehat{V}(y) \, \widehat{u_0}(p-y)}{p-y} \right] \right] \right\|_{L^\infty(\R \times [a,b])} \s^{-2} \; ;
		\end{align*}
		Note that the $L^\infty$-norms are well-defined because the functions are supported on the bounded domain $[p_1 + a, p_2 + b] \times [a, b]$ and continuous thanks to the hypotheses on $\widehat{V}$ and $\widehat{u_0}$. Since the $y$-integral in $W_{\infty, 2}^s(p)$ is actually defined over the bounded interval $[a,b]$, the last right-hand side is $(\s,y)$-integrable over $[s, +\infty) \times [a,b]$ and independent from $p$. By classical arguments on parametric integrals, we deduce that $W_{\infty, 2}^s$ is a $\mathcal{C}^1$-function on $\R$.
		
		By addition, we deduce that the function $W_\infty^s = W_{\infty, 1}^s + W_{\infty,2}^s$ is also continuously differentiable with a support contained in $[p_1 + a, p_2 + b]$.
		\item By the definition of $W_{\infty, 2}^s(p)$, we have
		\begin{align}
		 	\Big| W_{\infty, 2}^s(p)\Big|
		 		& \leqslant \frac{1}{8} \int_a^b \left| \partial_y \left[ \frac{1}{p-y} \partial_y \left[ \frac{1}{p-y} \partial_y \frac{\widehat{V}(y) \, \widehat{u_0}(p-y)}{p-y} \right] \right] \right| dy \int_s^\infty \s^{-3} ds \nonumber \\
		 		& \leqslant \frac{1}{16} \int_a^b \left| \partial_y \left[ \frac{1}{p-y} \partial_y \left[ \frac{1}{p-y} \partial_y \frac{\widehat{V}(y) \, \widehat{u_0}(p-y)}{p-y} \right] \right] \right| dy \, s^{-2} \; . \label{eq:w2_ineq_0}
		\end{align}
		Using similar arguments to those of the preceding point, we prove that
		\begin{equation} \label{eq:def_r}
			r(p) := \frac{1}{16} \int_a^b \left| \partial_y \left[ \frac{1}{p-y} \partial_y \left[ \frac{1}{p-y} \partial_y \frac{\widehat{V}(y) \, \widehat{u_0}(p-y)}{p-y} \right] \right] \right| dy
		\end{equation}
		defines a continuous function with compact support, so it is in particular integrable, and verifies inequality \eqref{eq:w2_ineq} according to \eqref{eq:w2_ineq_0}.
	\end{enumerate}	
\end{proof}

In the following proposition, we prove that $W_\infty^s(p)$ defined in Lemma \ref{lem:prop_winfty} is indeed the limit of $W(t,p)$ as $t$ tends to infinity. This is achieved by providing an upper bound for the absolute value of the difference between these two terms, the upper bound tending to $0$ as $t$ tends to infinity. Note that the regularity of the Fourier transforms of the potential and of the initial datum is again used in the following proof to carry out integrations by parts.
	
\begin{proposition} \label{prop:exp_winfty}
	Let $k \geqslant 3$ and $\ell \geqslant 4$ be two integers and suppose that $u_0 \in H^3(\R)$ and $V \in L^2(\R)$ satisfy respectively Conditions \emph{($\mathcal{I}_{[p_1,p_2]}^\ell$)} and \emph{($\mathcal{P}_{[a,b]}^k$)}. Assume in addition that $0 \notin [p_1, p_2]$.\\
	Let $t > 0$ and $s \in (0, t)$. Then the function $W(t,.)$ introduced in Proposition \ref{prop:rewriting} verifies 
	\begin{equation*}
		\forall \, p \in \R \qquad \big| W(t,p) - W_\infty^s(p) \big| \leqslant r(p) \, t^{-2} \; ,
	\end{equation*}
	where the functions $W_\infty^s$ and $r$  have been introduced in Lemma \ref{lem:prop_winfty} .
\end{proposition}

\begin{proof}
	Let $t > s > 0$. We first split the $\s$-integral in $W_{\infty, 2}^s(p)$ as follows:
	\begin{equation*}
		W_{\infty, 2}^s(p) = \frac{1}{8} \int_s^t \int_a^b \dots \, + \frac{1}{8} \int_t^\infty \int_a^b \dots \; ,
	\end{equation*}
	and we integrate then three times by parts the first term of the preceding sum to obtain:
	\begin{align*}
		& \frac{1}{8} \int_s^t \int_a^b \partial_y \left[ \frac{1}{p-y} \partial_y \left[ \frac{1}{p-y} \partial_y \frac{\widehat{V}(y) \, \widehat{u_0}(p-y)}{p-y} \right] \right] e^{-i \s (p-y)^2} dy \, \s^{-3} e^{i \s p^2} d\s \\
		& \hspace{5mm} = -i \int_s^t \int_a^b \widehat{V}(y) \, \widehat{u_0}(p-y) \, e^{-i \s (p-y)^2} \, dy \, e^{i \s p^2} d\s \; .
	\end{align*}
	It follows then
	\begin{align*}
		W_\infty^s(p) 
			& = W_{\infty, 1}^s(p) + W_{\infty, 2}^s(p) \\
			& = -i \int_0^s \int_a^b \widehat{V}(y) \, \widehat{u_0}(p-y) \, e^{-i \s (p-y)^2} \, dy \, e^{i \s p^2} d\s \\
			& \hspace{.5cm} - i \int_s^t \int_a^b \widehat{V}(y) \, \widehat{u_0}(p-y) \, e^{-i \s (p-y)^2} \, dy \, e^{i \s p^2} d\s \\
			& \hspace{1cm} + \frac{1}{8} \int_t^\infty \int_a^b \partial_y \left[ \frac{1}{p-y} \partial_y \left[ \frac{1}{p-y} \partial_y \frac{\widehat{V}(y) \, \widehat{u_0}(p-y)}{p-y} \right] \right] e^{-i \s (p-y)^2} dy \, \s^{-3} e^{i \s p^2} d\s \\
			& = -i \int_0^t \int_a^b \widehat{V}(y) \, \widehat{u_0}(p-y) \, e^{-i \s (p-y)^2} \, dy \, e^{i \s p^2} d\s \\
			& \hspace{.5cm} + \frac{1}{8} \int_t^\infty \int_a^b \partial_y \left[ \frac{1}{p-y} \partial_y \left[ \frac{1}{p-y} \partial_y \frac{\widehat{V}(y) \, \widehat{u_0}(p-y)}{p-y} \right] \right] e^{-i \s (p-y)^2} dy \, \s^{-3} e^{i \s p^2} d\s \\
			& = W(t,p) + W_{\infty,2}^t(p) \; .
	\end{align*}
	We apply finally Lemma \ref{lem:prop_winfty} to obtain
	\begin{equation*}
		\big| W(t,p) - W_\infty^s(p) \big| \leqslant \big| W_{\infty,2}^t(p) \big| \leqslant r(p) \, t^{-2} \; .
	\end{equation*}
\end{proof}

\begin{remark}
	The term $W_\infty^s(p)$ being the limit of $W(t,p)$ as $t$ tends to infinity, it does not depend in particular on $s > 0$.  So it will be denoted by $W_\infty(p)$ in the rest of the paper.
\end{remark}

\subsection{Time-asymptotic expansion for the second term of the Dyson-Phillips series}

In this last subsection, we exploit the results from the two preceding subsections to derive a time-asymptotic expansion of $S_2(t) u_0$ reflecting its spatial propagation. We first show that this term is time-asymptotically close to a free wave packet. In particular, the amplitude of this wave packet is actually given by the above limit $W_\infty$. Expanding this free wave packet as in Corollary \ref{cor:s1_exp} provides the desired expansion for the term $S_2(t) u_0$.\\

We start by proving that the term $S_2(t) u_0$ is time-asymptotically close to a free wave packet. To do so, we use the inequality given in Proposition \ref{prop:exp_winfty}.

\begin{proposition} \label{prop:first_exp_s2}
	Let $k \geqslant 3$ and $\ell \geqslant 4$ be two integers and suppose that $u_0 \in H^3(\R)$ and $V \in L^2(\R)$ satisfy respectively Conditions \emph{($\mathcal{I}_{[p_1,p_2]}^\ell$)} and \emph{($\mathcal{P}_{[a,b]}^k$)}. Assume in addition that $0 \notin [p_1, p_2]$.\\
	Then we have for all $(t,x) \in (0, +\infty) \times\R$,
	\begin{equation*}
		\left| S_2(t) u_0(x) - \frac{1}{2 \pi} \int_{\R} W_\infty(p) \, e^{-itp^2 + ixp} \, dp \right| \leqslant \frac{1}{2 \pi} \, \big\| r \big\|_{L^1} \, t^{-2} \; ,
	\end{equation*}
	where the functions $W_\infty$ and $r$ have been introduced in Lemma \ref{lem:prop_winfty}.
\end{proposition}

\begin{proof}
	Let $(t,x) \in (0, +\infty) \times\R$. Combining Propositions \ref{prop:rewriting} and \ref{prop:exp_winfty} leads to
	\begin{align*}
		& \left| S_2(t) u_0(x) - \frac{1}{2 \pi} \int_{\R} W_\infty(p) \, e^{-itp^2 + ixp} \, dp \right| \\
		& \hspace{1cm} = \left| \frac{1}{2 \pi} \int_{\R} W(t,p) \, e^{-itp^2 + ixp} \, dp - \frac{1}{2 \pi} \int_{\R} W_\infty(p) \, e^{-itp^2 + ixp} \, dp \right| \\
		& \hspace{1cm} \leqslant \frac{1}{2 \pi} \int_{\R} \big| W(t,p) - W_\infty(p) \big| \, dp \\
		& \hspace{1cm} \leqslant \frac{1}{2 \pi} \int_\R r(p) \, dp \, t^{-2} \; .
	\end{align*}
	Note that the last integral is finite since $r \in L^1(\R)$ according to Lemma \eqref{lem:prop_winfty}.
\end{proof}

We are now in position to obtain a time-asymptotic expansion of the term $S_2(t) u_0$. The proof of the following result consists mainly in expanding to one term the free wave packet introduced in Proposition \ref{prop:first_exp_s2}. This is achieved by applying Theorem \ref{thm:oi_exp} whose hypotheses are satisfied since the amplitude $W_\infty$ has been proved to be continuously dif\-fe\-ren\-tiable function on $\R$ with compact support. The resulting expansion is then put into the inequality of Proposition \ref{prop:first_exp_s2} to obtain at the end a time-asymptotic expansion for $S_2(t) u_0$.

\begin{theorem} \label{thm:s2_exp}
	Let $p_1$, $p_2$, $\widetilde{p}_1$ and $\widetilde{p}_2$ be four finite real numbers such that $[p_1, p_2] \subset \big(\widetilde{p}_1, \widetilde{p}_2 \big)$.
	Let $k \geqslant 3$ and $\ell \geqslant 4$ be two integers and suppose that $u_0 \in H^3(\R)$ and $V \in L^2(\R)$ satisfy respectively Conditions \emph{($\mathcal{I}_{[p_1,p_2]}^\ell$)} and \emph{($\mathcal{P}_{[a,b]}^k$)}. Assume in addition that $0 \notin [p_1, p_2]$.
	Consider the function $W_\infty$ defined in Lemma \ref{lem:prop_winfty} and define
	\begin{align*}
		& \bullet \quad J_{W_\infty}(t,x) = \frac{1}{2 \pi} \int_\R W_\infty(p) \, e^{-itp^2 + ixp} \, dp \qquad \forall \, (t,x) \in \R \times \R \; ; \\
		& \bullet \quad t_2^* = \argmin_{\tau \in \R} \left( \int_\R x^2 \, \big| J_{W_\infty}(t,x) \big|^2 \, dx - \frac{1}{\| W_\infty \|_{L^2(\R)}^2} \Big( \int_\R x \, \big| J_{W_\infty}(t,x) \big|^2 \, dx \Big)^2 \right) ; \\
		& \bullet \quad x_2^* = \frac{1}{\| W_\infty \|_{L^2(\R)}^2} \int_\R x \, \Big| J_{W_\infty}(t_2^*,x) \Big|^2 \, dx \; .
	\end{align*}
	Then for all $\displaystyle (t,x) \in \mathfrak{C}\big( [\widetilde{p}_1 + a, \widetilde{p}_2 + b], (t_2^*, x_2^*) \big)$ with $t \neq 0$, we have
	\begin{align*}
		& \left| S_2(t)u_0(x) - \frac{1}{\sqrt{2 \pi}} \, e^{- sgn(t-t_2^*) i \frac{\pi}{4}} \, e^{-it \big(\frac{x-x_2^*}{t - t_2^*}\big)^2 + ix \frac{x-x_2^*}{t - t_2^*}} \, W_\infty \left( \frac{x-x_2^*}{t - t_2^*} \right) |t-t_2^*|^{-\frac{1}{2}} \right| \\
		 & \hspace{5mm} \leqslant C_1(\delta, \tilde{p}_1 + a, \tilde{p}_2 + b) \, \sqrt{\int_\R x^2 \, \big| J_{W_\infty}(t_2^*,x) \big|^2 \, dx - \frac{1}{\| W_\infty \|_{L^2(\R)}^2} \Big(\int_\R x \, \big| J_{W_\infty}(t_2^*,x) \big|^2 \, dx \Big)^2} \, |t-t_2^*|^{-\delta} \\
		 & \hspace{1cm} + \frac{1}{2 \pi} \, \big\| r \big\|_{L^1} \, t^{-2} \; ,
	\end{align*}
	where the real number $\delta$ is arbitrarily chosen in $\big( \frac{1}{2}, \frac{3}{4} \big)$. And for all $\displaystyle (t,x) \in \mathfrak{C}\big( [\widetilde{p}_1 + a, \widetilde{p}_2 + b], (t_2^*, x_2^*) \big)^c$ with $t \neq 0$, we have
	\begin{align*}
		& \Big| S_2(t)u_0(x) \Big| \\
		& \hspace{1cm} \leqslant \Bigg( C_2(p_1, p_2, \tilde{p}_1 + a, \tilde{p}_2 + b) \, \sqrt{\int_\R x^2 \, \big| J_{W_\infty}(t_2^*,x) \big|^2 \, dx - \frac{1}{\| W_\infty \|_{L^2(\R)}^2} \Big( \int_\R x \, \big| J_{W_\infty}(t_2^*,x) \big|^2 \, dx \Big)^2} \\
		& \hspace{2cm} +  C_3(p_1, p_2, \tilde{p}_1 + a, \tilde{p}_2 + b) \, \big\| W_\infty \big\|_{L^\infty(\R)} \Bigg) \, |t-t_2^*|^{-1} + \frac{1}{2 \pi} \, \big\| r \big\|_{L^1} \, t^{-2} \; .
	\end{align*}
	All the above constants are explicitly given in \citet[Thm. 1.1]{dewez2020} and the integrable function $r$ is defined in Lemma \ref{lem:prop_winfty}.
\end{theorem}

\begin{proof}
	We define first the following free wave packet:
	\begin{equation*}
		 J_{\widetilde{W}_\infty}(t,x) := \frac{J_{W_\infty}(t,x)}{\| W_\infty \|_2} = \frac{1}{2 \pi} \int_\R \widetilde{W}_\infty(p) \, e^{-itp^2 + ixp} \, dp \; ,
	\end{equation*}
	where $\widetilde{W}_\infty := \frac{W_\infty}{\| W_\infty \|_2}$. The amplitude $\widetilde{W}_\infty$ is a $L^2$-normalised and continuously differentiable function with support contained in $[p_1 + a, p_2 + b]$ according to Lemma \ref{lem:prop_winfty}. The hypotheses of Theorem \ref{thm:oi_exp} are then verified and we obtain the following time-asymptotic expansion of the wave packet $J_{\widetilde{W}_\infty}$ in the cone $\mathfrak{C}\big( [\widetilde{p}_1 + a, \widetilde{p}_2 + b], (t_2^*, x_2^*) \big)$, 
	\begin{align*}
		& \left| J_{\widetilde{W}_\infty}(t,x) - \frac{1}{\sqrt{2 \pi}} \, e^{- sgn(t-t_2^*) i \frac{\pi}{4}} \, e^{-it \big(\frac{x-x_2^*}{t - t_2^*}\big)^2 + ix \frac{x-x_2^*}{t - t_2^*}} \, \widetilde{W}_\infty \left( \frac{x-x_2^*}{t - t_2^*} \right) |t-t_2^*|^{-\frac{1}{2}} \right| \\
		 & \hspace{5mm} \leqslant C_1(\delta, \tilde{p}_1 + a, \tilde{p}_2 + b) \, \sqrt{\int_\R x^2 \, \Big| J_{\widetilde{W}_\infty}(t_2^*,x) \Big|^2 \, dx -  \Big(\int_\R x \, \Big| J_{\widetilde{W}_\infty}(t_2^*,x) \Big|^2 \, dx \Big)^2} \, |t-t_2^*|^{-\delta} \; ,
	\end{align*}
	and the following uniform estimate outside,
	\begin{align*}
		\Big| J_{\widetilde{W}_\infty}(t,x) \Big|
			& \leqslant \Bigg( C_2(p_1, p_2, \tilde{p}_1 + a, \tilde{p}_2 + b) \, \sqrt{\int_\R x^2 \, \Big| J_{\widetilde{W}_\infty}(t_2^*,x) \Big|^2 \, dx - \Big( \int_\R x \, \Big| J_{\widetilde{W}_\infty}(t_2^*,x) \Big|^2 \, dx \Big)^2} \\
			& \hspace{1.5cm} +  C_3(p_1, p_2, \tilde{p}_1 + a, \tilde{p}_2 + b) \, \Big\| \widetilde{W}_\infty \Big\|_{L^\infty(\R)} \Bigg) \, |t-t_2^*|^{-1} \; .
	\end{align*}
	Combining finally the two preceding inequalities multiplied by $\| W_\infty \|_2$ and the inequality given in Proposition \ref{prop:first_exp_s2} via the triangle inequality provides the desired results.
\end{proof}

\appendix
\section{Appendix: Generic results from classical and functional analyses} \label{sec:appendix}

In this last section, we provide some generic results from classical and functional analyses which have been used in the present paper.\\

We start with the following lemma on which the proof of Proposition \ref{prop:admissible_potentials} is based. Lemma \ref{lem:localisation_fourier} provides a family of functions localised in a bounded frequency band $[a,b]$ and, at the same time, approximately localised in space in an interval centred on a point $x_0 \in \R$ with arbitrary precision if the band is sufficiently large. This result, whose proof lies essentially on Chebyshev's inequality, is originated from \cite{am2013}.

\begin{lemma} \label{lem:localisation_fourier}
	Let $k \geqslant 1$ be an integer, let $a$, $b$ and $x_0$ be three finite real numbers such that $a < b$, and let $\varphi$ be a $\mathcal{C}^k$-function such that $\supp \, \varphi \subseteq [-1,1]$.
	Let $f$ be the element of $L^2(\R)$ whose Fourier transform $\widehat{f}$ is the complex-valued function given by
	\begin{equation*}
		\forall \, p \in \R \qquad \widehat{f}(p) := \varphi \hspace{-1mm} \left( \frac{2p - (a+b)}{b-a} \right) e^{-i x_0 p} \; .
	\end{equation*}
	Then $\widehat{f}$ is a  $\mathcal{C}^k$-function on $\R$ supported on the interval $[a,b]$ and $f$ is an analytic function on $\R$ satisfying
	\begin{equation} \label{eq:chebyshev}
		\forall \, c > 0 \qquad \int_{|x-x_0| \geqslant c} \big| f(x) \big|^2 dx \leqslant \frac{2}{c^2} \, \frac{1}{b-a} \, \big\| \varphi' \big\|_{L^2(\R)}^2 \; .
	\end{equation}
\end{lemma}

\begin{proof}
	Since $\varphi: \R \longrightarrow \C$ is a  $\mathcal{C}^k$-function on $\R$ supported on $[-1,1]$, the Fourier transform $\widehat{f}$ of $f$ is clearly a $\mathcal{C}^k$-function on $\R$ such that
	\begin{equation*}
		\supp \, \widehat{f} \subseteq [a,b] \; ;
	\end{equation*}
	the boundedness of the support of $\widehat{f}$ implies in particular that $f$ is analytic on $\R$.\\
	Now let us prove inequality \eqref{eq:chebyshev}. For this purpose, we apply Chebyshev's inequality to the function $f$:
	\begin{equation} \label{eq:tcheby}
		\int_{|x-x_0| \geqslant c} \big| f(x) \big|^2 dx \leqslant \frac{1}{c^2} \, \int_{\R}  (x-x_0)^2 \, \big| f(x) \big|^2 dx  = \frac{1}{c^2} \, \int_{\R} (x-x_0)^2 \, \Big| \Big( \tf_{p \rightarrow x}^{-1} \widehat{f} \Big) (x) \Big|^2 dx \; ,
	\end{equation}
	for all $c > 0$. Then by a simple substitution, we have for all $x \in \R$,
	\begin{equation*}
		\Big( \tf_{p \rightarrow x}^{-1} \widehat{f} \Big)(x) = \frac{1}{2 \pi} \int_{\R} \varphi \hspace{-1mm} \left( \frac{2p - (a+b)}{b-a} \right) e^{i(x-x_0)p} \, dp  = \frac{b-a}{4 \pi} \, e^{i \frac{a+b}{2} (x-x_0)} \, \widehat{\varphi} \hspace{-1mm} \left( \frac{b-a}{2} \, (x_0-x) \right) \; .
	\end{equation*}
	Putting this into inequality \eqref{eq:tcheby} provides finally for all $c > 0$,
	\begin{align*}
		\int_{|x-x_0| \geqslant c} \big| f(x) \big|^2 dx
			& \leqslant \frac{1}{c^2} \, \frac{(b-a)^2}{16 \pi^2} \, \int_{\R} (x-x_0)^2 \left| \widehat{\varphi} \left( \frac{b-a}{2} \, (x_0-x) \right) \right|^2 dx \\
			& = \frac{1}{c^2} \, \frac{1}{2 \pi^2 (b-a)} \, \int_{\R} \big| y \, \widehat{\varphi}(y) \big|^2 dy \\
			& = \frac{2}{c^2} \, \frac{1}{b-a} \, \big\| \varphi' \big\|_{L^2(\R)}^2 \; ;
	\end{align*}
	note that we have used the substitution $y = \frac{b-a}{2} \, (x_0-x)$ to obtain the first equality and the classical relation $\widehat{(\varphi')}(y)= i \, y \, \widehat{\varphi}(y)$ together with Plancherel's theorem to obtain the second one.
\end{proof}

In the rest of the present section, we recall some results from semigroup theory and functional analysis which are used in Section \ref{sec:dyson}.\\
Let us remark that the following results are not proved but quoted from the li\-te\-ra\-tu\-re containing their proofs. Furthermore the operators $A,B,C$ and the semigroups $\big( T(t) \big)_{t \geqslant 0}$, $\big( S(t) \big)_{t \geqslant 0}$ used here are ge\-ne\-ric and do not refer to the particular objects which are considered in the preceding sections of this paper.\\

We recall first the notion of a classical solution for an abstract evolution equation; see \citet[Chap. II, Def. 6.1]{en2000}. Further we recall that if an operator generates a semigroup on a Banach space, then the classical solution of the evolution equation given by this operator exists, is unique and corresponds to the orbit of the initial value under the semigroup; see \citet[Chap. II, Prop. 6.2]{en2000}.

\begin{defprop} \label{defprop:existence}
	Consider the initial value problem
	\begin{equation} \label{eq:evol-eq-2}
		\left\{ \begin{array}{rl}
				& \hspace{-2mm} \dot{u}(t) = A \, u(t) \\ [2mm]
				& \hspace{-2mm} u(0) = v
		\end{array} \right. \; ,
	\end{equation}
	for $t \geqslant 0$, where $A : D(A) \subset X \longrightarrow X$ is the generator of a semigroup $\big( T(t) \big)_{t \geqslant 0}$ on the Banach space $X$.\\
	A function $u : [0, +\infty) \longrightarrow X$ is called a classical solution of \eqref{eq:evol-eq-2} if $u$ is continuously differentiable with respect to $X$, $u(t) \in D(A)$ for all $t \geqslant 0$, and $u$ satisfies \eqref{eq:evol-eq-2}.\\	
	If $v \in D(A)$ then the function
	\begin{equation*}
		u : t \in [0, +\infty) \longmapsto u(t) = T(t) \, v \; ,
	\end{equation*}
	is the unique classical solution of \eqref{eq:evol-eq-2}.
\end{defprop}

In the following theorem, we recall that the sum of a generator of a semigroup and a bounded operator on a Banach space generates a semigroup as well; see \citet[Chap. III, Thm. 1.3]{en2000}.

\begin{theorem}  \label{thm:perturb}
	Let $\big( A, D(A) \big)$ be the generator of a strongly continuous semigroup on a Banach space $X$. If $B$ is a bounded operator from $X$ into itself, i.e. $B \in \mathcal{L}(X)$, then the operator $\big( C , D(C) \big) := \big( A + B , D(A) \big)$ generates a strongly continuous semigroup on $X$.
\end{theorem}

The Dyson-Phillips series is now introduced in a generic setting: it provides a representation as a series of the semigroup generated by the operator $\big( A + B , D(A) \big)$; see \citet[Chap. III, Thm. 1.10]{en2000}.

\begin{theorem}  \label{thm:dyson}
	Let $\big( A, D(A) \big)$ be the generator of a strongly continuous semigroup $\big( T(t) \big)_{t \geqslant 0}$ on a Banach space $X$, and let $B \in \mathcal{L}(X)$. The strongly continuous semigroup $\big( S(t) \big)_{t \geqslant 0}$ generated by $\big( C , D(C) \big) := \big( A + B , D(A) \big)$ can be obtained as
	\begin{equation*} \label{eq:dyson-series}
		\lim_{N \rightarrow + \infty} \left\| \, S(t) - \sum_{n = 1}^{N} S_n(t) \, \right\|_{\mathcal{L}(X)} = \, 0 \; ,
	\end{equation*}
	where $S_1(t) := T(t)$ and
	\begin{equation*} \label{eq:dyson-integral}
		\forall \, v \in X \qquad S_{n+1}(t) v := \int_0^t S_n(t-\s) \, B \, T(\s) v \, d\s \; .
	\end{equation*}
\end{theorem}

In the final result, we recall that Bochner-type integration and the application of bounded operators can be interchanged; see \citet[Chap. V, Sec. 5, Cor. 2]{yosida1980}. This result has been used in the present paper to evaluate the terms of the Dyson-Phillips at any point $x \in \R$.

\begin{proposition} \label{prop:bochner}
	Let $A$ be a bounded operator acting between two Banach spaces $X$ and $Y$ and let $J \subseteq \R$ be an interval. If $F : J \longrightarrow X$ is a Bochner-integrable function, then $A F : J \longrightarrow Y$ is also a Bochner-integrable function and
	\begin{equation*}
		A \left( \int_J F(s) \, ds \right) = \int_J A \, F(s) \, ds \; .
	\end{equation*}
\end{proposition}

\bibliography{biblio}

\end{document}